\documentclass[12pt]{article}
\usepackage{amsmath,amssymb,theorem}
\usepackage{graphicx}
\usepackage{caption}
\usepackage{subfigure}
\usepackage{psfrag}
\usepackage{epstopdf}
\newtheorem{thm}{Theorem}[section]
\newtheorem{conj}[thm]{Conjecture}
\newtheorem{lem}[thm]{Lemma}
\theorembodyfont{\rmfamily}
\def\pf{\bigskip\noindent {\bf Proof.}~~}

\def\dfn#1{{\sl #1}}
\def\mytextindent#1{\indent\llap{#1\enspace}\ignorespaces}

\def\myitem{\par\hangindent\parindent\mytextindent}

\newcommand{\less}{\backslash}

\newcounter{counter}
\def\pr#1{(\ref{#1})}

\begin{document}
\title{Coloring graphs with forbidden minors }
\author{Martin Rolek\thanks{E-mail address: mrolek@knights.ucf.edu.}
~and~ Zi-Xia Song\thanks{Corresponding author. E-mail address: Zixia.Song@ucf.edu.}\\
Department of Mathematics\\
University of Central Florida\\
 Orlando, FL 32816
}

\maketitle
\begin{abstract}

 Hadwiger's conjecture from 1943 states that for every integer $t\ge1$,  every graph either can be $t$-colored or has a subgraph that can be contracted to the complete graph on $t+1$ vertices.  
 As pointed out by Paul Seymour in his recent survey on Hadwiger's conjecture,   proving that graphs with no $K_7$ minor are $6$-colorable is the first case of Hadwiger's  conjecture that is still open.
It  is not  known yet whether  graphs with no $K_7$ minor are $7$-colorable.
Using a Kempe-chain argument  along with the fact that an induced path on three vertices is dominating in a graph with independence number two,  we first give a very short and computer-free proof of a recent result of  Albar and Gon\c calves and generalize it to the next step by showing  that every graph  with no $K_t$ minor is  $(2t-6)$-colorable, where $t\in\{7,8,9\}$.  We then prove that  graphs with no $K_8^-$ minor are $9$-colorable and  graphs with no $K_8^=$ minor are $8$-colorable. Finally we prove that if Mader's bound for the extremal function for $K_p$ minors is true, then every graph with no $K_p$ minor is $(2t-6)$-colorable  for all $p\ge5$.  This  implies our first result. We believe that the Kempe-chain method we have developed in this paper is of independent interest.

\end{abstract}

\baselineskip=18pt
\section{Introduction}

All graphs in this paper are finite and simple. A graph  $H$ is a \dfn{minor} of a graph $G$ if  $H$ can be
 obtained from a subgraph of $G$ by contracting edges.  We write $G>H$ if 
$H$ is a minor of $G$.
In those circumstances we also say that  $G$ has an $H$ \dfn{minor}.  
\medskip

Our work is motivated by the following Hadwiger's conjecture~\cite{hc}, which is perhaps the most famous conjecture in graph theory, as pointed out by Paul Seymour in his recent survey~\cite{Seymour}.

\begin{conj}\label{hc}  For every integer $t\ge1$, every graph with no 
$K_{t+1}$ minor is $t$-colorable.
\end{conj}

Hadwiger's conjecture is trivially true for $t\le2$, and reasonably easy for
$t=3$, as shown by Dirac~\cite{Dirac1952}. However, for $t\ge4$, Hadwiger's conjecture
implies the Four Color Theorem. 
Wagner~\cite{wagner} proved that the case $t=4$ of Hadwiger's conjecture is, in fact,
 equivalent to the Four Color Theorem, and the same was shown for $t=5$
by Robertson, Seymour and Thomas~\cite{RST}.   Hadwiger's conjecture remains open for $t\ge6$.   As pointed out by Paul Seymour~\cite{Seymour} in his recent survey on Hadwiger's conjecture, proving that graphs with no $K_7$ minor are $6$-colourable is thus the first case of Hadwiger's  conjecture that is still open.
It  is not even known yet whether  every graph with no $K_7$ minor is $7$-colorable. 
 Kawarabayashi and Toft~\cite{ktoft} proved that every graph with no $K_7$ or $K_{4,\, 4}$ minor is $6$-colorable. 
 Jakobsen~\cite{Jakobsen1972, Jakobsen1983} proved that every graph with no $K_7^{=}$ minor is  $6$-colorable and  every graph with no $K_7^{-}$ minor is  $7$-colorable,  where for any integer $p>0$,
 $K_p^{-}$ (resp. $K_p^{=}$) denotes the graph obtained from $K_p$ by removing one edge (resp. two edges).    For more information on Hadwiger's conjecture, the readers are referred to an earlier  survey   by Toft~\cite{Toft}   and a very recent  informative survey  due to Seymour~\cite{Seymour}.  \medskip

   Albar and Gon\c calves~\cite{AG2015} recently proved the following:

\begin{thm}\label{K7K8} (Albar and Gon\c calves~\cite{AG2015}) 
Every graph with no $K_7$ minor is $8$-colorable, and every graph with no $K_8$ minor is $10$-colorable
\end{thm}

The proof of Theorem~\ref{K7K8} is computer-assisted and not simple. In this  paper,  we  apply a Kempe-chain argument (see Lemma~\ref{wonderful} below) along with the fact that an induced path on three vertices is dominating in a graph with independence number two in order to  give a much shorter and computer-free proof of Theorem~\ref{K7K8}.  In addition, we generalize it to the next step by proving the following. 

\begin{thm}\label{main}
Every graph with no $K_t$ minor is $(2t-6)$-colorable, where  $t\in\{7, 8,9\}$.
\end{thm}

We want to point out  that our proof of Theorem~\ref{main}  does not rely on  Mader's deep result on the connectivity of contraction-critical graphs (see Theorem~\ref{7con} below). Theorem~\ref{main} states that (i) every graph with no $K_7$ minor is $8$-colorable;  (ii) every graph with no $K_8$ minor is $10$-colorable; and (iii)  every graph with no $K_{9}$ minor is $12$-colorable. We prove Theorem~\ref{main} in Section~\ref{Kt}.  \medskip

Applying the method we developed in the proof of Theorem~\ref{main} and Mader's deep result (Theorem~\ref{7con}), we  then prove  two new results  Theorem~\ref{mainK8-} and Theorem~\ref{mainK8=}.

\begin{thm}\label{mainK8-}
Every graph with no  $K_8^-$  minor is  $9$-colorable.
\end{thm}

 \begin{thm}\label{mainK8=}
Every graph with no $K_8^=$ minor is $8$-colorable.
\end{thm}
 
Our proofs of  Theorem~\ref{mainK8-}  and Theorem~\ref{mainK8=}  are both short and computer-free and will be presented in Section~\ref{K8-} and  Section~\ref{K8=}, respectively. \medskip
 
 To end this paper, we first propose a conjecture in Section~\ref{remarks}.  We then apply Lemma~\ref{wonderful}  to prove that if  Conjecture~\ref{conj3} ( see Section~\ref{remarks}) is true, 
 then every graph with no $K_p$ minor is $(2t-6)$-colorable  for all $p\ge5$.   Our proof  of the last result does not rely on the connectivity of contraction-critical graphs and the new idea  we introduced   yields  a different/short  proof of Theorem~\ref{main}.\medskip

 To prove our results, we need to investigate the basic properties of contraction-critical graphs. For a positive integer $t$, a graph $G$ is \dfn{$t$-contraction-critical} if $\chi(G)=t$ and any proper minor of $G$ is $(t-1)$-colorable.   
 Lemma~\ref{dirac} below  is  a folklore result 
 which is an extension of Dirac's initial work~\cite{dirac2} on contraction-critical graphs. A proof of  Lemma~\ref{dirac} can be easily obtained from the definition of $k$-contraction-critical graphs. 

\begin{lem}\label{dirac}(Dirac~\cite{dirac2})   Every   $k$-contraction-critical graph $G$ satisfies the following:

\myitem {(i)} for any $v\in V(G)$, $\alpha(N(v))\le d(v)-k+2$, where $\alpha(N(v))$ denotes the independence nunber of the subgraph of $G$ induced by $N(v)$;  
\myitem {(ii)} no minimal separating set of $G$ is  a clique.
 \end{lem}

\begin{figure}[htb]
\centering
\includegraphics[width=220px]{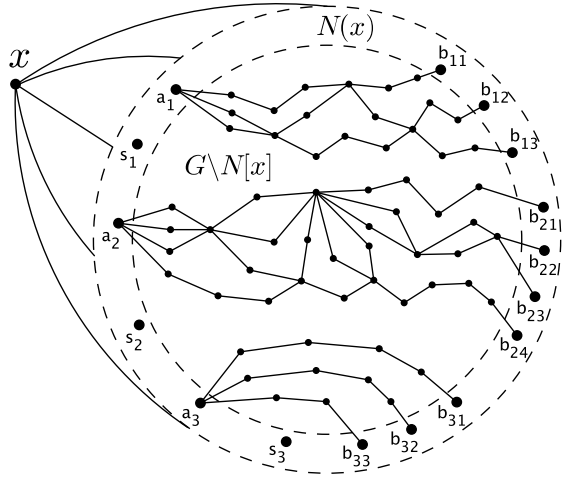}
\caption{An application of Lemma~\ref{wonderful} to $x$ with $S = \{s_1, s_2, s_3\}$ and $M = \{\{a_1 b_{11}, a_1 b_{12}, a_1 b_{13}\}, \{a_2 b_{21}, a_2 b_{22}, a_2b_{23}, a_2b_{24}\}, \{a_3 b_{31}, a_3b_{32}, a_3b_{33}\}\}$.}
\label{wonderfullemma}
\end{figure}

 Lemma~\ref{wonderful} below on contraction-critical graphs  turns out to be very powerful, as the existence of pairwise vertex-disjoint paths  is guaranteed without using the connectivity of such  graphs.  
  If two vertices  $u,v$ in a graph $G$  are not adjacent,  we say that $uv$  is  a \dfn{missing edge} of $G$.  One possible  application of Lemma~\ref{wonderful} is depicted in Figure~\ref{wonderfullemma}.

 \begin{lem}\label{wonderful}  Let $G$ be   any $k$-contraction-critical graph. Let   $ x\in V(G)$ be a vertex  of  degree $k+s$ with $\alpha(N(x))=s+2$ and let $S\subset N(x)$ with $|S|=s+2$ be any independent set, where  $k\ge4$  and $s\ge0$ are integers.  If $N(x)\less S$ is not a clique, then  for any  $M=\{\{a_1b_{11}, \dots, a_1b_{1r_1}\}, \{a_2b_{21}, \dots, a_2b_{2r_2}\}, \dots, \{a_mb_{m1}, \dots, a_mb_{mr_m}\}\}$, where $m, r_i\ge1$,  $r_1+r_2+\dots+r_m+ m\le k-2$, the vertices $a_1, \dots, a_m,
  b_{11}, \dots, b_{mr_m}\in N(x)\less S$ are all distinct, and  for any $i\in \{1,2,\dots, m\}$, the set $\{a_ib_{i1}, \dots, a_ib_{ir_i}\}$ consists of $r_i$ missing edges of $N(x)\less S$ with $a_i$ as a common end,   then for each  $i\in \{1,2,\dots, m\}$ there exist paths $P_{i1}, \dots, P_{ir_i}$  in $G$ such that each $P_{ij}$ has ends $a_i, b_{ij}$ and  all its internal vertices in $G\less N[x]$ for all $j=1,2,\dots, r_i$.  Moreover, for any $1\le i< \ell \le m$, the paths $P_{i1}, \dots, P_{ir_i}$ are vertex-disjoint from the paths $P_{\ell1}, \dots, P_{\ell r_{\ell}}$.  \end{lem} 

\pf  Let $G$, $x$, $S$ and $M$ be as given in the statement.  Let $H$ be obtained from $G$ by contracting $S\cup\{x\}$ into a single vertex, say  $w$.
Then $H$ is $(k-1)$-colorable. Let $c : V(H)\rightarrow \{1,2, \dots, k-1\}$ be a proper  $(k-1)$-coloring of $H$. We may assume that $c(w)=1$.
Then  each of the colors $ 2,  \dots, k-1$ must appear in $N(x)\less S$, else we could assign $x$ the missing color and all vertices in $S$ the color $1$ to  obtain a proper $(k-1)$-coloring of $G$, a contradiction.
Since $|N(x)\less S|=k-2$,   we have $c(u)\ne c(v)$ for any two distinct vertices $u, v$ in $N(x)\less S$.   
 We next claim that  for each $i\in\{1,2, \dots, m\}$ and each $j\in\{1,2,\dots, r_i\}$   there must exist a path between $a_i$ and $b_{ij}$ with its internal vertices in $G\less N[x]$. Suppose not. Let $i\in\{1,2, \dots, m\}$ and $j\in\{1,2, \dots, r_i\}$  be such that there is no such path between $a_i$ and $b_{ij}$.
Let $H^*$ be the subgraph of $H$ induced by the vertices colored  $c(a_i)$ or $c(b_{ij})$ under the coloring $c$. Then $V(H^*)\cap N(x)=\{a_i,b_{ij}\}$. Notice that  $a_i$ and $b_{ij}$ must belong to different components of $H^*$ as there is no  path between $a_i$ and $b_{ij}$ with its internal vertices in $G\less N[x]$.
By switching the colors on the component of $H^*$ containing $a_i$, we obtain a $(k-1)$-coloring of $H$ with the color $c(a_i)$ missing on $N(x)\less S$, a contradiction.
This proves that  there must exist a path $P_{ij}$ in $H^*$ with ends $a_i, b_{ij}$ and all its internal vertices in $H^*\less N[x]$ for  each $i\in\{1,2, \dots, m\}$ and each $j\in\{1,2,\dots, r_i\}$.
Clearly, for any $1\le i< \ell \le m$, the paths $P_{i1}, \dots, P_{ir_i}$ are vertex-disjoint from the paths $P_{\ell1}, \dots, P_{\ell r_{\ell}}$, because  no two  vertices of $a_1, \dots, a_r,
  b_{11}, \dots, b_{mr_m}$ are colored the same under the coloring $c$.  \hfill\vrule height3pt width6pt depth2pt \\

\noindent {\bf Remark.}  If $r_1=r_2=\dots=r_m=1$ in the statement of Lemma~\ref{wonderful}, we simply write $M=\{a_1b_{11}, a_2b_{21}, \dots, a_mb_{m1}\}$, and so $M$ is a matching of missing edges of $N(x)\less S$.
In this case, the paths $P_{11}, P_{21}, \dots, P_{m1}$ are pairwise vertex-disjoint if $m\ge2$. If $m=1$ in the statement of Lemma~\ref{wonderful}, we simply write $M=\{a_1b_{11}, \dots, a_1b_{1r_1}\}$.
In this case, the paths $P_{11}, \dots, P_{1r_1}$ have $a_1$ as a common end and are not necessarily pairwise internally vertex-disjoint if $r_1\ge2$.  \medskip

Theorem~\ref{7con} below  is a  deep result of Mader~\cite{7con} and will be used in the proofs of Theorem~\ref{mainK8-} and Theorem~\ref{mainK8=}.   It seems very difficult to improve Theorem~\ref{7con} for small values of $k$.  For larger values of $k$, some better results can be found.  Kawarabayashi~\cite{Kawa2007} has shown that any minimal non-complete $k$-contraction-critical graph with no $K_k$ minor is $\left\lceil 2k/27\right\rceil$-connected, while Kawarabayashi and Yu~\cite{Kawa2013} have shown that any minimal such graph is $\left\lceil k/9 \right\rceil$-connected. Chen, Hu and Song~\cite{CHS} recently improved the bound further  by showing that any minimal such  graph is $\left\lceil k/6 \right\rceil$-connected. 

\begin{thm}\label{7con} (Mader~\cite{7con}) 
For $k \ge 7$, every $k$-contraction-critical graph is $7$-connected.
\end{thm}

We also need the following lemma  in the proofs of Theorem~\ref{mainK8-} and Theorem~\ref{mainK8=}. 

\begin{lem}\label{lem:7Conn2K6}
For any $7$-connected graph $G$, if $G$ contains two different $K_6$-subgraphs, then $G > K_8^-$.
\end{lem}

\pf  Let $H_1, H_2$ be two different $K_6$-subgraphs of $G$ with   $V(H_1) = \{v_1, \dots, v_6\}$ and $V(H_2) = \{w_1, \dots, w_6\}$.  Let $t=|V(H_1)\cap V(H_2)|$. Then $0\le t\le5$. We may assume that   $v_i = w_i$ for all $i\le t$ if $t\ne0$.
Assume that $t= 5$. Then $H_1 \cup H_2$ has a $K_7^-$-subgraph of $G$.  Since $G$ is $7$-connected, it is easy to see that $G > K_8^-$ by contracting a component of $G\less (H_1\cup H_2)$ into a single vertex.
So we may assume that $t \le 4$. 
Then there exist $6-t$  pairwise disjoint paths $P_{t+1}, \dots, P_{6}$  between $H_1\less H_2$ and $H_2\less H_1$ in $G\less (V(H_1)\cap V(H_2))$. We may assume that $P_i$ has ends $v_{i}, w_{i}$ for all $i=t+1,  \dots,6$.  
Then $G \less \{v_1, \dots, v_5, w_6\}$ is connected, so there must exist a path $Q$ with one end, say $u$, in $(P_{t+1}\less v_{t+1} )\cup \dots \cup (P_5\less v_5)$, the other end, say $v$,  in $P_6 \less w_6$, and no internal vertices in any of $\{v_1, \cdots, v_t\}, P_{t+1}, \dots, P_6$.
We may assume that $u$ lies on the path $P_5\less v_5$. 
Let $P_5^*$ be the subpath of $P_5$ with ends $u, w_5$, and $P_6^*$ be the subpath of $P_6$ with ends $v, v_6$.
Now contracting $P_5^*$ onto $w_5$, $P_5\less P_5^*$ onto $v_5$, $P_6^*$ and $Q\less u$ onto $v_6$, $P_6\less P_6^*$ onto $w_6$, and  each of $P_{t+1}, \dots,  P_4$ to a single vertex, together with $v_1, \dots, v_t$ if $t\ne0$, yields  a $K_8^-$ minor in $G$, as desired. \hfill\vrule height3pt width6pt depth2pt\\

We need to introduce more notation.   For a graph $G$ we use $|G|$ and $\delta (G)$ to denote the number
of vertices and minimum degree of $G$, respectively. For a subset $S$ of $V(G)$,
the subgraph induced by $S$ is denoted by
$G[S]$ and $G\less S =G[V(G)\less S]$.
The degree and neighborhood of a vertex $v$ in $G$ are denoted by 
 $d(v)$ and $N(v)$, respectively.  By abusing
notation we will also denote by $N(v)$ the graph induced by the set
$N(v)$.  We define $N[v]=N(v)\cup \{v\}$, and similarly will use the
same symbol for the graph induced by that set.    For $S\subseteq V(G)$, if $G[S]$ is connected, then we denote by $G/S$ the graph obtained from $G$ by contracting $G[S]$  into a single vertex and deleting all resulting parallel edges and loops.
For $A, B\subseteq V(G)$, we say that $A$ is \dfn{complete} to $B$ if each vertex in $A$ is adjacent to all vertices in $B$, and $A$ is \dfn{anti-complete} to $B$ if no vertex in $A$ is adjacent to any vertex in $B$.
If $A=\{a\}$, we simply say $a$ is complete to $B$ or $a$ is anti-complete to $B$.  The {\dfn{join}} $G+H$ (resp. {\dfn{union}} $G\cup H$) of two 
vertex disjoint graphs
$G$ and $H$ is the graph having vertex set $V(G)\cup V(H)$  and edge set $E(G)
\cup E(H)\cup \{xy\, |\,  x\in V(G),  y\in V(H)\}$ (resp. $E(G)\cup E(H)$). \medskip


\section{Proof of Theorem~\ref{main}:  coloring $K_t$-minor free graphs}\label{Kt}

Results on the extremal function for $K_t$ minors will be needed to prove Theorem~\ref{main}.   Before doing so, we need to define $(H_1, H_2, k)$-cockade. For  graphs $H_1, H_2$ and an integer
$k$, let us define
an {\it $(H_1, H_2,k)$-cockade} recursively as follows. Any graph isomorphic to $H_1$ or $H_2$
is an $(H_1, H_2,k)$-cockade. Now let $G_1$, $G_2$ be $(H_1, H_2,k)$-cockades and let $G$
be obtained from the disjoint union of $G_1$ and $G_2$ by identifying a clique
of size $k$ in $G_1$ with a clique of the same size in $G_2$. Then the graph
$G$ is also an $(H_1, H_2,k)$-cockade, and every $(H_1, H_2,k)$-cockade can be constructed in this way.
If $H_1=H_2=H$, then $G$ is simply called an $(H, k)$-cockade. 
   The following  Theorem~\ref{k7} was first shown by Dirac~\cite{Dirac1964} for $p \le 5$ and by Mader~\cite{mader} for $p = 6, 7$.

\begin{thm}\label{k7}(Mader~\cite{mader}) 
For every integer $p=1,2,\ldots,7$, a graph on $n\ge p$ vertices
and at least $(p-2)n-{p-1\choose2}+1$ edges has a $K_p$ minor.
\end{thm}

The edge bound in Theorem~\ref{k7} is referred to as Mader's bound for the extremal function for $K_p$ minors.  J\o rgensen~\cite{Jorgensen1994} and later  the second author and Thomas~\cite{SongThomas2006} generalized Theorem~\ref{k7} to $p=8$ and $p=9$, respectively,  as follows.  

\begin{thm}\label{k8}(J\o rgensen~\cite{Jorgensen1994})  
Every graph on $n \ge 8$ vertices with at least $6n - 20$ edges either has a  $K_8$  minor or is a $(K_{2,2,2,2,2}, 5)$-cockade.
\end{thm}

\begin{thm}\label{k9}(Song  and Thomas~\cite{SongThomas2006}) 
Every graph on $n \ge 9$ vertices with at least $7n - 27$ edges either has a  $K_9$  minor, or is  a $(K_{1,2,2,2,2,2}, 6)$-cockade, or is isomorphic to $K_{2,2,2,3,3}$.
\end{thm}

It seems hard to generalize Theorem~\ref{k7} for all values of $p$.   In 2003, Seymour and Thomas~\cite{SongThomas2006} proposed the following conjecture. 

\begin{conj}\label{conj2}(Seymour and Thomas~\cite{SongThomas2006}) For every $p\ge1$ there exists a constant $N=N(p)$ such that
every $(p-2)$-connected graph on $n\ge N$ vertices and at least
$(p-2)n-{p-1\choose2}+1$ edges has a $K_p$ minor.
\end{conj} 

By Theorem~\ref{k9}, Conjecture~\ref{conj2} is true for $p\le9$.  \medskip

We next prove the following Lemma~\ref{k_tk1}, which can be obtained  from   the (computer-assisted) proof of Lemma 3.7 in~\cite{SongThomas2006}. Here we give a computer-free proof of Lemma~\ref{k_tk1} so that  the  proof of Theorem~\ref{main} is  also computer-free.

 \begin{lem}\label{k_tk1} For $7\le t\le 9$, let $H$ be a graph with $2t-5$ vertices and  
$\alpha(H)=2$. Then $H>K_{t-2}\cup K_1$. \end{lem} 

\setcounter{counter}{0}

\pf  Suppose that $H$ has no $K_{t-2}\cup K_1$ minor.   Then $\omega(H)\le t-3$. We claim that  \medskip

\noindent \refstepcounter{counter}\label{e:noK_{t-3}} (\arabic{counter})  
$\omega(H)\le t-4$.  

\pf Suppose that  $\omega(H)= t-3$. Let $K\subseteq H$ be isomorphic to $K_{t-3}$. Then $|H\less K|=t-2\ge5$. If $H\less K$ contains an induced  $3$-path,  say $P$,  with ends $y,z$, where a $3$-path is a path with three vertices,  then every vertex of $K$ is adjacent to either $y$ or $z$ because $\alpha(H)=2$. By contracting the path $P$ into a single vertex, we see that $H[K\cup P]>K_{t-2}$ and so $H>K_{t-2}\cup K_1$, a contradiction. Thus $H\less K$ does not contain an induced path on three vertices.  Since $\alpha(H)=2$, it follows that  $H\less K$ is a disjoint union of two cliques, say $A_1$ and $A_2$. For $i=1,2$, let $K_i=\{v\in V(K): v \text{ is not adjacent to some vertex in  } A_{3-i}\}$. Since $\alpha(H)=2$, $K_i$ is complete to $A_i$ for each $i$. Thus $H\less (K_i\cup A_i)$ is a clique for each $i\in\{1,2\}$ and so either  $H\less (K_1\cup A_1)$ or $H\less (K_2\cup A_2)$ is a clique of size at least $t-2$, contrary to the fact that $\omega(H)\le t-3$. \hfill\vrule height3pt width6pt depth2pt \medskip

 Let $q=\delta(H)$ and let $y\in V(H)$ be  a vertex with $d(y)=q$. Let $J=H\less N[y]$. Since $\alpha(H)=2$,   $J$ is a clique of size $2t-q-6$.  By \pr{e:noK_{t-3}}, $|J|=2t-q-6\le t-4$ and so  $q\ge t-2$. We next show that  \medskip

\noindent \refstepcounter{counter}\label{e:2cliques} (\arabic{counter})  
for any $A\subseteq N(y)$ with $|A|\ge6$,  either $H[A\cup\{y\}]$ contains two vertex-disjoint $3$-paths or $H[A]$ is a disjoint union of two cliques. 

\pf Suppose  $H[A]$ is  not a disjoint union of two cliques. Then  $H[A]$ is connected because $\alpha(H)=2$.  We next show that $H[A\cup\{y\}]$  contains two vertex-disjoint induced $3$-paths.   By  \pr{e:noK_{t-3}},    $H[A]$ is not a clique and thus  contains an induced $3$-path, say $abc$, with ends $a,c$. Let $\{d_1, d_2, \dots, d_s\}=A\less\{a,b,c\}$, where $s=|A|-3\ge3$.  Clearly $H[A\cup\{y\}]$  contains two vertex-disjoint induced $3$-paths if $H[\{d_1, d_2, \dots, d_s\}]$ is not a clique, since $yd_i$ is an edge for $1 \le i \le s$.
So we may assume that $H[\{d_1, d_2, \dots, d_s\}]=K_s$. First assume that $a$ is complete to $ \{d_1, d_2, \dots, d_s\}$. 
By  \pr{e:noK_{t-3}}, $b$ is not complete to $\{d_1, d_2, \dots, d_s\}$. We may assume that $bd_1\notin E(H)$. Clearly $H[\{a,y,c\}]$ and $H[\{d_1, b, d_i\}]$ are two vertex-disjoint $3$-paths if $bd_i\in E(H)$ for some $i\ne 1$.
So we may assume that $bd_i\notin E(H)$. Now either $H[\{b, a, d_1\}]$ and $H[\{c, y, d_2\}]$ (if $cd_2\notin E(H)$) or $H[\{ a, d_2, c\}]$ and $H[\{b, y, d_1\}]$ (if $cd_2\in E(H)$) are  two vertex-disjoint $3$-paths. Next assume that $a$ is not complete to $ \{d_1, d_2, \dots, d_s\}$.  We may assume that $ad_1\notin E(H)$.
Then $cd_1\in E(H)$ because $\alpha(H)=2$. By symmetry, we may assume that $cd_2\notin E(H)$.
Then $ad_2\in E(H)$.
Now either $H[\{c, d_1, d_2\}]$ and $H[\{a, y, d_3\}]$ (if $ad_3\notin E(H)$) or $H[\{ a, d_3, d_1\}]$ and $H[\{c, y, d_2\}]$ (if $ad_3\in E(H)$) are  two vertex-disjoint $3$-paths, as desired. \hfill\vrule height3pt width6pt depth2pt \medskip

\noindent \refstepcounter{counter}\label{e:q=t-2} (\arabic{counter})  
$q=t-2$.  

\pf  Suppose  $q\ge t-1$. By Theorem~\ref{k7},  $(t-4)(2t-6)-{t-3\choose 2}\ge e(H\less y)\ge q|H|/2-q=q(|H|-2)/2\ge(t-1)(2t-7)/2$, which yields that  $t=9$ and  $q=t-1=8$.    
Then $H$ is a graph on  thirteen vertices. 
 Clearly,  $J=K_4$.  Let $z\in N(y)$ be such that $|N(z)\cap V(J)|$ is maximum. Since    $e_H(J, N(y))\ge20$, we have $|N(z)\cap V(J)|\ge3$.  If $|N(z)\cap V(J)|=4$, then $H[\{z\}\cup V(J)]=K_5$ and $|N(y)\less z|=7$.
Clearly $H>K_7\cup K_1$ if $N[y]\less z$ has two vertex-disjoint induced $3$-paths.
By   \pr{e:2cliques}, $N[y]\less z$ is thus a disjoint union of two cliques, say $A_1, A_2$.
By   \pr{e:noK_{t-3}}, we may assume that $A_1=K_3$ and $A_2=K_4$.
Let $a\in A_1$.
By   \pr{e:noK_{t-3}} again, $a$ is not complete to $\{z\} \cup J$   and thus $d_H(a)\le 7$,  contrary to the fact that $q=8$. 
 Thus $|N(z)\cap V(J)|=3$.
Let $z'\in V(J)$ be the non-neighbor of $z$.  By the choice of $z$,  every vertex in $N(y)$ has at least one non-neighbor in $J$ and so $\delta(N(y))\ge 4$. Since $d(z)\ge8$, $|N(z)\cap N(y)|\ge 4$. By \pr{e:noK_{t-3}}, $N(z)\cap N(y)$ is not a clique and so $z'$ is adjacent to at least one vertex, say $w$, in $N(z)\cap N(y)$, because $\alpha(H)=2$. Now the edge $zw$ is dominating $J$, i.e., every vertex in $J$ is adjacent to either $z$ or $w$.
Notice that $|N(y)\less\{z,w\}|=6$. If  $N[y]\less\{z,w\}$ contains two vertex-disjoint induced $3$-paths, say $P_1$ and $P_2$, then $H>K_7\cup K_1$ by contracting the edge $zw$ and the two $3$-paths $P_1$ and $P_2$ into three distinct vertices, respectively,  a contradiction.
Thus $N[y]\less\{z,w\}$ does not contain two vertex-disjoint induced $3$-paths.
By   \pr{e:2cliques},  $N(y)\less\{z,w\}$ is a disjoint union of two cliques, say $B_1$ and $B_2$. Since $\delta(N(y))\ge4$, we have $B_1=B_2=K_3$. By  \pr{e:noK_{t-3}},   $H[B_1\cup\{z,w, y\}]$ is not a clique.
Let $w'\in B_1$ be such that either $ww'\notin E(H)$ or $zw'\notin E(H)$.
Since $w'$ is adjacent to at most three vertices of $J$, we see that $d_H(w')\le 7$, contrary to the fact that $q=8$. 
 \hfill\vrule height3pt width6pt depth2pt \medskip

By  \pr{e:q=t-2},  $q=t-2$. If $t=7$, then $H$ is a graph on nine vertices with $\delta(H)=5$. Thus there exists a vertex $z\in V(H)$ such that $d_H(z)\ge6$ and so $N[z]$ contains a $K_4$-subgraph because $\alpha(N(z))=2$, contrary to  \pr{e:noK_{t-3}}.
Hence $t\ge8$.  Now $H\less N[y]$ is a clique of size $t-4$ and $|N(y)|=t-2\ge6$.  Clearly  $H>K_{t-2}\cup K_1$ if  $H\less N(y)$ contains two vertex-disjoint  induced 3-paths,   a contradiction.
Thus by   \pr{e:2cliques},  $N(y)$ is a disjoint union of two cliques, say $A_1$ and $A_2$. For $i=1,2$, let $K_i=\{v\in H\less N[y]: v \text{ is not adjacent to some vertex in  } A_{3-i}\}$. Since $\alpha(H)=2$, $K_i$ is complete to $A_i$ for each $i$.
Thus $H\less (K_i\cup A_i\cup\{y\})$ is a clique for each $i\in\{1,2\}$ and so at least one of them is of size at least $t-3$, contrary to   \pr{e:noK_{t-3}}.
This completes the proof of Lemma~\ref{k_tk1}.  \hfill\vrule height3pt width6pt depth2pt\\

We are now ready to prove Theorem~\ref{main}. \bigskip

\setcounter{counter}{0}

\noindent {\bf Proof of Theorem~\ref{main}.}  \,Suppose the assertion is false. Let $G$ be a  graph with no $K_t$ minor such  that $G$ is  not $(2t-6)$-colorable.  
We may choose such a  graph $G$ so that it is $(2t-5)$-contraction-critical.
Let $x\in V(G)$ be of minimum degree.  Since $K_{2,2,2,3,3}$ and each  $(K_{2,2,2,2,2}, 5)$-cockade are $5$-colorable, and every $(K_{1, 2,2,2,2,2}, 6)$-cockade is $6$-colorable, it follows from Theorem~\ref{k7}, Theorem~\ref{k8} and Theorem~\ref{k9} that  $d(x)\le 2t-5$.  On the other hand, since $G$ is $(2t-5)$-contraction-critical, by Lemma~\ref{dirac}(i) applied to $N(x)$, we see that $d(x)\ge(2t-5)-2+\alpha(N(x))$.  Clearly,  $\alpha(N(x))\ge2$, otherwise $N[x]$  is a clique of size at least $2t-5>t$, contrary to the fact that   $G$  has no $K_t$ minor. Thus  $\alpha(N(x))\ge2$ and so $d(x)\ge2t-5$.   Then  $d(x)=2t-5\ge t+2$. By Lemma~\ref{dirac}(i) applied to $N(x)$ again, we have $\alpha(N(x))=2$. We next show that  \medskip

\noindent \refstepcounter{counter}\label{e:noK_{t-1}} (\arabic{counter})  
$G$ has no $K_{t-1}$-subgraph. 

\pf  Suppose $G$ contains $K_{t-1}$ as a subgraph. Let $H\subseteq G$ be isomorphic to  $K_{t-1}$. Since $\delta(G)=d(x)\ge t+2$, every vertex in $H$ is adjacent to at least one vertex in $G\less H$. Then $G\less H$ is disconnected, since otherwise $G>K_t$ by contracting $G\less H$ into a single vertex, a contradiction. Let $G_1$ be a component of $G\less H$. Then $N(G_1) :=\{v\in V(H):  v \, \text{ is adjacent to a vertex in}\, G_1\}$ is a minimal separating set of $G$.  In particular,  $N(G_1)$ is a clique, contrary to Lemma~\ref{dirac}(ii).  \hfill\vrule height3pt width6pt depth2pt \medskip

\noindent \refstepcounter{counter}\label{e:cnbr} (\arabic{counter})  
for any $u\in N(x)$, $|N(x)\cap N(u)|\ge t-3$.

\pf  Suppose that there exists a vertex   $u\in N(x)$ such that  $|N(x)\cap N(u)|\le t-4$. Since $\alpha(N(x))=2$,  $N(x)$ contains a clique of size  $|N(x)\backslash N[u]|\ge t-2$ and so  $N[x]$ has a $K_{t-1}$-subgraph, contrary to \pr{e:noK_{t-1}}.  \hfill\vrule height3pt width6pt depth2pt \medskip

By Lemma~\ref{k_tk1}, $N(x)>K_{t-2}\cup K_1$. Let $y\in N(x)$ be  such that $N(x)\backslash y>K_{t-2}$. Clearly, $y$ is not adjacent to every vertex in $N(x)\less y$, otherwise $G>N[x]>K_t$, a contradiction. Let  $\{y_1, \dots, y_p\}=N(x)\less N[y]$, where $p=2t-5-|N(x)\cap N[y]|$. Then $y$ is not adjacent to $y_1, y_2, \dots, y_p$. By \pr{e:noK_{t-1}} and   \pr{e:cnbr},    $N[y]\cap N(x)$ is not a clique.  Let $uw$ be a missing edge in $N(y)$. By Lemma~\ref{wonderful} applied  to $N(x)$  with $k=2t-5$,  $S=\{u,w\}$ and  $M=\{yy_1, yy_2, \dots, yy_p\}$,    there exists a path $P_i$  with ends  $y$ and $y_i$ and all  its internal vertices in $G\less N[x]$ for each $i\in\{1,2,\dots, p\}$. 
  Note that the paths $P_{1}\less y_1, \dots, P_{p}\less y_p$ have $y$ as a common end. By contracting all $P_{i}\less y_i$ onto $y$, we see that $G>K_{t}$, a contradiction. \hfill\vrule height3pt width6pt depth2pt \bigskip
 
 
\section{Proof of Theorem~\ref{mainK8-}: coloring $K_8^-$-minor free graphs}\label{K8-}

The primary purpose of this section is to prove Theorem~\ref{mainK8-} which states that every graph with no $K_8^-$ minor  is $9$-colorable. We need the following results. Theorem~\ref{thm:K567-Extremal} was proved by Dirac~\cite{Dirac1964} for the cases $p = 5, 6$, and by Jakobsen~\cite{Jakobsen1983} for the case $p = 7$.

\begin{thm}\label{thm:K567-Extremal}(Dirac~\cite{Dirac1964}, Jakobsen~\cite{Jakobsen1972,Jakobsen1983}) 
For $p=5, 6, 7$, if $G$ is a graph with $n \ge p$ vertices and at least $(p - \frac{5}{2})n - \frac{1}{2}(p - 3)(p - 1)$ edges, then $G > K_p^-$, or $G$ is a $(K_{p - 1}, p - 3)$-cockade when $p \ne 7$, or $p = 7$ and $G$ is a $(K_{2,\, 2,\, 2,\, 2}, \, K_6, \, 4)$-cockade.
\end{thm} 

 Jakobsen~\cite{Jakobsen1983} also conjectured that Theorem~\ref{thm:K567-Extremal} extends to $p = 8$, which was confirmed by the second author~\cite{Song2005}: 

\begin{thm}\label{thm:K8-Extremal}(Song~\cite{Song2005}) 
If $G$ is a graph with $n \ge 8$ vertices and at least $\frac{1}{2}(11n - 35)$ edges, then $G > K_8^-$ or $G$ is a $(K_{1,\, 2,\, 2,\, 2,\, 2},\,  K_7,\, 5)$-cockade. 
\end{thm} 

The extremal function for $K_p^-$ minors remains open for $p\ge9$.  The following lemma from \cite{Song2005} will also be needed. 

\begin{figure}[htb]
\centering
\includegraphics[width=150px]{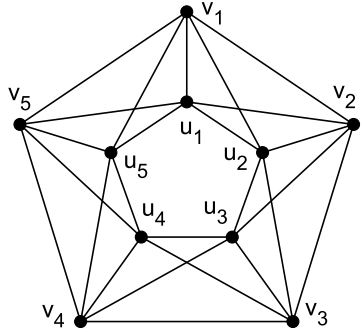}
\caption{The graph $J$.}
\label{fig:J}
\end{figure}

\begin{lem}\label{lem:delta5Nbrhd}(Song~\cite{Song2005}) 
Let $G$ be a graph with $8 \le |G| \le 10$ and $\delta(G) \ge 5$. Then either $G > K_6^- \cup K_1$ or $G$ is isomorphic to one of $\overline{C_8}$, $\overline{C_4} + \overline{C_4}$, $\overline{K_3} + C_5$, $\overline{K_2} + \overline{C_6}$, $K_{2,\, 3,\, 3}$, or $J$, where $J$ is the graph depicted in Figure \ref{fig:J}. In particular, all of these graphs are edge maximal (subject to not having a $K_6^- \cup K_1$-minor) with maximum degree $\le |G| - 2$. Moreover, $\overline{C_8} > K_6$, $\overline{C_4} + \overline{C_4} > K_6$, and $J > K_6$.
\end{lem}

Notice that of the counterexamples listed in Lemma~\ref{lem:delta5Nbrhd}, only the graph $J$ has ten vertices, and none has exactly nine vertices.  We first prove the following lemma.

\begin{lem}\label{lem:10vertexalpha2}
Let $G$ be a graph with $n = 10$ vertices and $\alpha(G) = 2$.  Then either $G > K_6^- \cup K_1$, or $G$ contains a $K_5 \cup K_5$  subgraph, or $G$ is isomorphic to the graph $J$ depicted in Figure \ref{fig:J}. \end{lem}

\pf If $\delta(G)\ge 5$, then by Lemma~\ref{lem:delta5Nbrhd}, either $G > K_6^- \cup K_1$ or $G$ is isomorphic to $J$.  So we may assume that $\delta(G)\le4$. 
 Let $x \in V(G)$ be such that $d(x)=\delta(G)$.  Since $\alpha(G) = 2$, one can easily see that  $G > K_6 \cup K_1$ if $d(x)\le3$.  Hence we may further assume that $d(x)=4$.  Then  $G\less N[x]$ must be a $K_5$ as  $\alpha(G) = 2$.  If $N[x]$ also induces a $K_5$, then $G$ contains a $K_5 \cup K_5$ subgraph.  Otherwise, some edge is missing from $N(x)$, say $y, z \in N(x)$ with $yz \notin E(G)$.  Then since $\alpha(G) = 2$, each vertex in $G\less N[x]$ must be adjacent to either  $y$ or $z$.  Thus  $G/\{x,y,z\} > K_6 \cup K_1$, as desired. This completes the proof of Lemma~\ref{lem:10vertexalpha2}.  \hfill\vrule height3pt width6pt depth2pt\bigskip

We are now ready to  prove Theorem~\ref{mainK8-}.

\setcounter{counter}{0}

 \pf   Let $G$ be a graph with no  $K_8^-$  minor.
Suppose for a contradiction that $\chi(G) \ge 10$.
We may choose such a  graph $G$ so that it is $10$-contraction-critical.  Let $x\in V(G)$ be of minimum degree. 
Since $G$ is $10$-contraction-critical and has  no $K_8^-$ minor, by Lemma~\ref{dirac}(i) applied to $N(x)$, we see that $\alpha(N(x))\ge2$ and $\delta(G)\ge10$.
On the other hand, since every $(K_{1,2,2,2,2}, K_7, 5)$-cockade is $7$-colorable, by Theorem~\ref{thm:K8-Extremal} we see that $\delta(G) \le 10$.   Thus $\delta(G)=10$.
By Lemma~\ref{dirac}(i)  applied to $N(x)$, we have \medskip

\noindent \refstepcounter{counter}\label{e:alpha} (\arabic{counter}) 
 $\alpha(N(x)) = 2$.\medskip

We next show that \medskip

\noindent \refstepcounter{counter}\label{e:NoJ} (\arabic{counter}) 
 $N(x)$ is not isomorphic to the graph $J$.

\pf Suppose that  $N(x)$ is isomorphic to the graph $J$.  Let  the vertices of $J$ be labeled as depicted in Figure \ref{fig:J}. By Lemma~\ref{wonderful} applied to $J$ with $S=\{v_2, v_5\}$ and $M=\{ \{u_1u_3, u_1u_4, u_1v_3, u_1v_4\}, \{u_2u_5\}\}$ with  $m=2$, $r_1=4, r_2=1$, there exist   paths $P_{11}$, $P_{12}$, $P_{13}$, $P_{14}$, $P_{21}$   such that the paths $P_{11}$, $P_{12}$, $P_{13}$, $P_{14}$, $P_{21}$   have  ends  $\{u_1, u_3\}$, $\{u_1, u_4\}$, $\{u_1, v_3\}$,  $\{u_1, v_4\}$, and $\{u_2, u_5\}$,   respectively, and all their internal vertices in $G\backslash N[x]$.  Moreover,  the paths $P_{11}$, $P_{12}$, $P_{13}$, $P_{14}$  are vertex-disjoint from the path $P_{21}$. By contracting $(P_{11}\less u_3)\cup  (P_{12}\less u_4)\cup(P_{13}\less v_3)\cup  (P_{14}\less v_4)$ onto $u_1$,   $P_{21}\less u_2$ onto $u_5$, and  $J[\{v_2, v_1, v_5\}]$ into a single vertex, we see that $G>K_8$, a contradiction.  
 \hfill\vrule height3pt width6pt depth2pt\medskip

\noindent \refstepcounter{counter}\label{e:K5UK5} (\arabic{counter})   
 $N(x)$ contains $K_5 \cup K_5$ as a subgraph.

\pf Suppose that $N(x)$ does not contain $K_5 \cup K_5$ as a subgraph.  Then by \pr{e:alpha}, \pr{e:NoJ} and Lemma~\ref{lem:10vertexalpha2}, we see that $N(x) > K_6^- \cup K_1$.
Let $y \in N(x)$ be a vertex such that $N(x) \setminus \{y\} > K_6^-$.  Clearly, $y$ is not adjacent to every vertex in $N(x)\less y$, otherwise $G>N[x]>K_8^-$, a contradiction. Let  $\{y_1, \dots, y_p\}=N(x)\less N[y]$, where $p=10-|N(x)\cap N[y]|\ge1$. Then $y$ is not adjacent to $y_1, y_2, \dots, y_p$. Clearly,   $N(x)\less \{y, y_i\}$ is not a clique for all $i\in\{1,2,\dots, p\}$. By Lemma~\ref{wonderful} applied $p$ times to $N(x)$  with $k=10$,  $s=0$ and $m=1$  (where $M=\{yy_i\}$ for $i=1,2\dots, p$),    there exists a path $P_i$  between $y$ and $y_i$ with its internal vertices in $G\less N[x]$ for each $i\in\{1,2,\dots, p\}$. 
 Note that the paths $P_{1}, \dots, P_{p}$ have $y$ as a common end. By contracting each $P_{i}\less y_i$ onto $y$, we see that $G>K_{8}^-$, a contradiction. \hfill\vrule height3pt width6pt depth2pt 
\medskip

 By \pr{e:K5UK5}, $x$ belongs to two different $K_6$-subgraphs of $G$.  By Theorem~\ref{7con}, $G$ is $7$-connected.  By Lemma~\ref{lem:7Conn2K6}, $G > K_8^-$.
This contradiction completes the proof of Theorem~\ref{mainK8-}.  \hfill\vrule height3pt width6pt depth2pt\\


 \section{Proof of Theorem~\ref{mainK8=}: coloring $K_8^=$-minor free graphs}\label{K8=}
 
We prove Theorem~\ref{mainK8=} in this section. The following result will be needed. Theorem~\ref{thm:K8=Extremal} for the cases  $p=5,6$ is due to Dirac~\cite{Dirac1964}, and Theorem~\ref{thm:K8=Extremal} for the cases $p=7,8$ is due to Jakobsen~\cite{Jakobsen, Jakobsen1972}.

\begin{thm}\label{thm:K8=Extremal}(Dirac~\cite{Dirac1964}, Jakobsen~\cite{Jakobsen, Jakobsen1972})  
For integer $p$ with  $5\le p\le8$,  every  graph  with $n \ge p$ vertices and at least $(p-3)n-\frac12(p-1)(p-4)$ edges either contains a $K_p^=$-minor or  is a $(K_{p-1}, p-4)$-cockade.
\end{thm}

We are ready to prove  Theorem~\ref{mainK8=}. 
\setcounter{counter}{0}

\pf  Suppose the assertion is false. Let $G$ be a  graph with no $K_8^=$ minor such  that $\chi(G) \ge 9$.  
We may choose such a  graph $G$ so that it is $9$-contraction-critical.
Let $x\in V(G)$ be of minimum degree. Since $G$ is $9$-contraction-critical and has no $K_8^=$ minor, by Lemma~\ref{dirac}(i) applied to $N(x)$, we see that $\alpha(N(x))\ge2$ and  $d(x)\ge9$.  On the other hand, since  each  $(K_{7}, 4)$-cockade is $4$-colorable, it follows from Theorem~\ref{thm:K8=Extremal} for $p=8$  that  $d(x)\le 9$. Thus $d(x)=9$,  and so $\delta(G) = 9$.   
It follows from Theorem~\ref{thm:K8=Extremal} for $p=8$ again  that \medskip

\noindent \refstepcounter{counter}\label{e:K8=:n9} (\arabic{counter})  
$G$ contains at least 28 vertices of degree 9.\medskip

Since $G$ has  no $K_8^=$ minor, by Lemma~\ref{dirac}(i) applied to $N(x)$,\medskip

\noindent \refstepcounter{counter}\label{e:K8=:alpha} (\arabic{counter})  
$\alpha(N(x)) = 2$.\medskip

We next show that \medskip

\noindent \refstepcounter{counter}\label{e:K8=:deltaNbr4} (\arabic{counter}) 
$N(x)$  contains  $K_5$ as a subgraph.

\pf  Suppose that  $N(x)$ does not contain $K_5$ as a subgraph. Then $\omega(N(x))\le4$ and by \pr{e:K8=:alpha}, $\delta(N(x))\ge4$. We claim  that $\delta(N(x))=4$. Suppose that $\delta(N(x))\ge5$. By Lemma~\ref{lem:delta5Nbrhd} applied to $N(x)$, we see that $N(x)>K_6^-\cup K_1$. Let $y \in N(x)$  be such that  $N(x)-y>K_6^-$.   Clearly $y$ has at least two non-neighbors in $N(x)-y$, otherwise $N[x]>K_8^=$, a contradiction. Let $y_1, y_2, \dots, y_j\in N(x)-y$ be all non-neighbors of $y$, where $j=|N(x)\backslash N[y]|\ge2$.
Since $\omega(N(x))\le4$, $N(x)\cap N(y)$ must have a missing edge, say $uv$. By Lemma~\ref{wonderful} applied   to $N(x)$ with $S=\{u, v\}$ and $M=\{yy_1, \dots, yy_j\}$, there exist $j$   paths $P_1, P_2, \dots, P_j$ such that each path $P_i$ has ends $\{y, y_i\}$ and all its internal vertices in $G\backslash N[x]$. By contracting all the edges of each $P_i\backslash y_i$ onto $y$ for all $i\in\{1,2,\dots, j\}$, we see that $G>K_8^-$, a contradiction. This proves that $\delta(N(x))=4$, as claimed. \medskip

Let $y \in N(x)$ be such that  $y$ has  degree four in $N(x)$ with  $e(N(y)\cap N(x))$ maximum.  Let  $Z=\{z_1, z_2, z_3,  z_4\}$ be the set of all neighbors of $y$ in $N(x)$.  Since $\omega(N(x))\le4$,  $N[y] \cap N(x)$ is not complete. We may assume that $z_1z_2\notin E(G)$. 
 By \pr{e:K8=:alpha},  $N(x)\less N[y]$ induce a $K_4$ subgraph.  Let $W=\{w_1, w_2, w_3, w_4\}=N(x)\less N[y]$.
We next show that \medskip

\noindent ~($*$)  each of $z_3, z_4$ has  at most one neighbor in $W$.

\pf Suppose, say $z_4$, is adjacent to at least two vertices in $W$. Then the subgraph induced on $W\cup \{z_4\}$ has a $K_5^=$ minor and thus  $N[x]>K_8^=$ if $z_3$ is adjacent to all vertices in $W$ (by contracting the path $z_1 y z_2$ into a single vertex), a contradiction. Thus we may assume that $z_3$ is not adjacent to $w_1, \dots, w_i$, where $1\le i\le 4$.  By Lemma~\ref{wonderful} applied  to $N(x)$ with $S=\{z_1, z_2\}$ and $M=\{z_3w_1, \dots, z_3w_i\}$, there exist $i$   paths $P_1, P_2, \dots, P_i$ such that for each $j=1,2,\dots, i$, the path $P_j$ has ends $\{z_3, w_j\}$ and all its internal vertices in $G\backslash N[x]$. By contracting  all $P_j\backslash w_j$ onto $z_3$  and the path $z_1yz_2$ into a single vertex, we see that $G>K_8^=$, a contradiction. This proves  ($*$).  \hfill\vrule height3pt width6pt depth2pt\medskip

We next claim that $N[y]\cap N(x)=K_5^-$. Suppose $z_3z_4\notin E(G)$. By symmetry, we may apply ($*$) to the missing edge $z_3z_4$ in $N(x)$, and so we see that each of $z_1,z_2$ has at most one neighbor in $W$.  Hence  $e_G(Z, W)\le 4$. On the other hand, since $\alpha(N(x))=2$, each $w_i$ must be adjacent to at least one of the vertices in $\{z_1, z_2\}$ and $\{z_3, z_4\}$, respectively, for all $i=1,2,3,4$. Thus $e_G(W, Z)\ge8$, a contradiction.
This proves that  $z_3z_4\in E(G)$ and thus $N(x)$ does not have two independent missing edges.  Next if $z_1z_3\notin E(G)$, then $z_2z_4, z_2z_3\in E(G)$ and $N[y]\cap N(x)=K_5^=$. Since $\omega(N(x))\le4$, we may assume that $z_1w_1\notin E(G)$. Then $w_1$ must be adjacent to both $z_2$ and $z_3$ by \pr{e:K8=:alpha}.
By  applying  ($*$) to the missing edges $z_1 z_2$ and $z_1 z_3$,  we see that  $\{z_2, z_3\}$ is anti-complete to $\{ w_2, w_3, w_4\}$ and $z_4$ has at most one neighbor in $W$.  By (2), $z_1$ is complete to  $ \{w_2, w_3, w_4\}$. 
  Since $z_4$ has at most one neighbor in $W$, we may assume that $w_4z_4\notin E(G)$. Now $w_4$ has degree four in $N(x)$ with  $N[w_4]\cap N(x)=K_5^-$, contrary to the choice of $y$. Thus $N[y]\cap N(x)=K_5^-$ with $z_1z_2$ the only missing edge, as claimed. 
\medskip

Since $\delta(N(x)) = 4$, each of $z_1$ and $z_2$ has at least one neighbor in $W$.
By \pr{e:K8=:alpha}, each of $w_1, \dots, w_4$ is adjacent to at  least one of $z_1, z_2$, and so either  $z_1$ or  $z_2$ has at least two neighbors in $W$.
By symmetry, we may assume that $z_1$ has more neighbors in $W$ than $z_2$.  On the other hand, each vertex in $Z$ has at least one non-neighbor in $W$ as $\omega(N(x)\le4$. Thus, $z_1$ has  either  two  or three neighbors in $W$. We  consider the following two cases.\medskip

\begin{figure}[htb]
\centering
\includegraphics[width=225px]{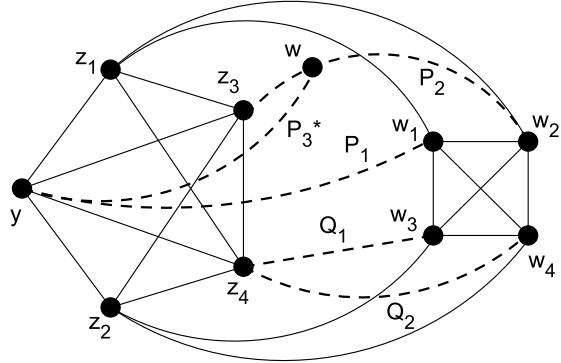}
\caption{$z_1$ has exactly two neighbors in $W$.}
\label{fig:case1}
\end{figure}

First, assume that $z_1$ has exactly two neighbors in $W$, say $w_1, w_2$.  Then $z_2$ must have exactly two neighbors   $w_3, w_4$ in $W$. By ($*$), each of $z_3, z_4$ has at most one neighbor in $W$. We may assume that $z_4$ is not adjacent to $w_3, w_4$, and $z_3w_2\notin E(G)$.  
By Lemma~\ref{wonderful} applied twice  to $N(x)$ with $S=\{z_1, z_2\}$ and $M\in \{\{yw_1, z_3w_2, z_4w_3\}, \{yw_2, z_4w_4\}\}$, there exist  three vertex-disjoint  paths $P_1, P_2, Q_1$ and two  vertex-disjoint paths $P_3, Q_2$  such that the paths $P_1$, $P_2$, $P_3$, $Q_1$, and  $Q_2$  have  ends $\{y, w_1\}$,  $\{z_3, w_2\}$,  $\{y, w_2\}$, $\{z_4, w_3\}$ and $\{z_4, w_4\}$,  respectively, and all their internal vertices in $G\backslash N[x]$, as depicted in Figure~\ref{fig:case1}.
Notice that each $P_i$ is vertex-disjoint from $Q_j$ for all $i\in\{1,2,3\}$ and $j\in\{1,2\}$, $P_1$ and $P_2$ are vertex-disjoint but $P_3$ and $P_2$ are  not necessarily  vertex-disjoint. If $P_3$ and $P_2$  have only $w_2$ in common, then contracting $P_1\less w_1, P_3\less w_2$ onto $y$, $P_2\less w_2$ onto $z_3$, and  $Q_1\less w_3, Q_2\less w_4$ onto $z_4$,  and $w_1w_3, w_2w_4$ into two distinct vertices  yields a $K_8^=$ minor in $G$, a contradiction.
Thus $P_3$ and $P_2$ must have an internal vertex in common. Let $w$ be the first vertex on $P_3$ (when $P_3$ is read from $y$ to $w_2$) that is also on $P_2$.  Then $w\notin V(P_1)$.  Let $P_3^*$ be the subpath of $P_3$ from $y$ to $w$, $P_2^*$ be the subpath of $P_2$ from $w$ to $w_2$.  Notice that $P_3^*\less w$ is vertex-disjoint from $P_2$ but not necessarily internally disjoint from $P_1$.
Now  contracting $P_1\less w_1$, $P^*_3\less w$ onto $y$, $P^*_2$ onto $w_2$, $P_2\less P_2^*$ onto $z_3$,  $Q_1\less w_3, Q_2\less w_4$ onto $z_4$, and $w_1w_3, w_2w_4$ into two distinct vertices yields another  $K_8^=$ minor in $G$, a contradiction. \medskip

\begin{figure}[htb]
\centering
\includegraphics[width=225px]{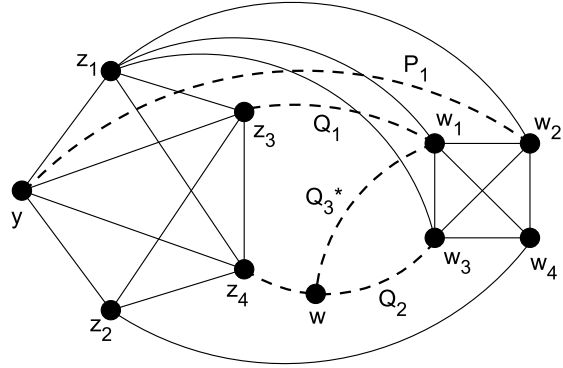}
\caption{$z_1$ has exactly three neighbors in $W$.}
\label{fig:case2}
\end{figure}

It remains to consider the case when  $z_1$ has exactly three neighbors, say $w_1, w_2, w_3$ in $W$.  
Then $z_2$ is adjacent to $w_4$. By ($*$), we may assume that $w_1$ is not adjacent to $z_3, z_4$, and $z_4w_3\notin E(G)$. 
By Lemma~\ref{wonderful} applied twice  to $N(x)$ with $S=\{z_1, z_2\}$ and $M\in \{\{yw_2, z_3w_1, z_4w_3\}, \{z_4w_1\}\}$, there exist  vertex-disjoint  paths $P_1, Q_1, Q_2$ and  another  path $Q_3$  such that the paths $P_1, Q_1, Q_2, Q_3$   have  ends $\{y, w_2\}$,  $\{z_3, w_1\}$, $\{z_4,w_3\}$, and $\{z_4, w_1\}$, respectively,  and all their internal vertices in $G\backslash N[x]$, as depicted in Figure~\ref{fig:case2}. Notice that $P_1$ is vertex-disjoint from $Q_j$ for $j = 1, 2, 3$,  but that $Q_3$ is not necessarily  internally vertex-disjoint from either $Q_1$ or $Q_2$. 
If $Q_3$ and $Q_2$ have only $z_4$ in common, then we obtain a $K_8^=$ minor by contracting $P_1$ and $z_2w_4$ into two distinct vertices,  $Q_1\less z_3$, $Q_3\less z_4$ onto $w_1$, and $Q_2\less z_4$ onto $w_3$, a contradiction.
Thus $Q_3$ and $Q_2$ must have an internal vertex in common. Let $w$ be the first vertex on $Q_3$ (when $Q_3$ is read from $w_1$ to $z_4$) that is also on $Q_2$. Then $w\notin V(Q_1)$.
Let $Q_3^*$ be the subpath of $Q_3$ from $w_1$ to $w$, $Q_2^*$ be the subpath of $Q_2$ from $w$ to $z_4$. Notice that $Q_3^*\backslash w$ is vertex-disjoint from $Q_2$ but not necessarily internally disjoint from $Q_1$.
Now we obtain another  $K_8^=$ minor by contracting $P_1$ and $z_2w_4$ into two distinct vertices,  $Q_1\less z_3$, $Q^*_3\less w$ onto $w_1$, $Q^*_2 $ onto $z_4$, and $Q_2\less Q_2^*$ onto $w_3$, a contradiction. This completes the proof of \pr{e:K8=:deltaNbr4}. \hfill\vrule height3pt width6pt depth2pt\medskip

By  \pr{e:K8=:deltaNbr4},  
every vertex of degree 9 belongs to some $K_6$-subgraph of $G$.
By \pr{e:K8=:n9}, $G$ contains at least five different $K_6$-subgraphs.  By Theorem~\ref{7con}, $G$ is $7$-connected and thus $G>K_8^-$ by Lemma~\ref{lem:7Conn2K6}. 
This contradiction completes the proof of Theorem~\ref{mainK8=}. \hfill\vrule height3pt width6pt depth2pt\\

 \section{Concluding Remarks}\label{remarks}

It seems very difficult  to prove that every graph with no $K_7$ minor is $7$-colorable. We establish in~\cite{RolekSong2} the properties of $8$-contraction-critical graphs with no $K_7$ minor  to shed some light on this open problem.  As pointed out by Robin Thomas (personal communication with the second author), one might be able to settle this open problem by using the key ideas in~\cite{RST}. However, it seems very hard to    prove that every  $8$-contraction-critical graph is $8$-connected, which means the truth of the Seymour-Thomas Conjecture (i.e. Conjecture~\ref{conj2}) can not be applied easily to solve Hadwiger's conjecture. Theorem~\ref{k7}  is such a nice result. We believe that  Mader's bound for the extremal function for $K_p$ minors is true as follows:  \medskip

\begin{conj}\label{conj3} For every $p\ge1$, 
every  graph $G$ on $n$ vertices and at least
$(p-2)n-{p-1\choose2}+1$ edges either has a $K_p$ minor or is  $(p-1)$-colorable.
\end{conj} 

By Theorem~\ref{k9}, Conjecture~\ref{conj3} is true for $p\le9$.  As mentioned earlier, Lemma~\ref{wonderful} turns out to be very powerful.  We believe that the application of Lemma~\ref{wonderful} that we have developed in this paper is of independent interest.  To end this section, we apply Lemma~\ref{wonderful} along with a new idea (namely, considering the chromatic number of $N(x)$)  to prove that the truth of Conjecture~\ref{conj3} implies that every graph with no $K_p$ minor is $(2p-6)$-colorable for all $p\ge5$.  Since  Conjecture~\ref{conj3} is true for $p\le9$, we see that Theorem~\ref{final}  implies Theorem~\ref{main}.  

\begin{thm}\label{final}
If Conjecture~\ref{conj3} is true, then every graph with no $K_p$ minor is $(2p-6)$-colorable for all $p\ge5$.
\end{thm}

\pf Suppose the assertion is false. 
By Theorem~\ref{k9}, we have $p\ge10$.
Among all minimum counterexamples, we choose $G$ so that $G$ has no $K_p$ minor and $G$ is $(2p - 5)$-contraction-critical. Let $x \in V(G)$ be such that $d(x) = \delta(G)$.  
By the assumed truth of Conjecture~\ref{conj3},  $d(x) \le 2p - 5$.  On the other hand, by Lemma~\ref{dirac}(i) applied to $N(x)$, we see that $\alpha(N(x))\ge2$ and  $d(x) \ge 2p - 5$.  Hence $d(x) = 2p - 5$.
By Lemma~\ref{dirac}(i) again, we have \medskip

\setcounter{counter}{0}

\noindent \refstepcounter{counter}\label{e:falpha} (\arabic{counter}) 
$\alpha(N(x)) = 2$. \medskip

Our strategy now will be to examine the subgraph $N(x)$ and its chromatic number.  We first show that  \vskip 0.1cm

\noindent \refstepcounter{counter}\label{e:fomega} (\arabic{counter}) 
$\omega(N(x)) \le p-3$, and so $\delta(N(x))\ge p-3$.

\pf  Suppose that  $\omega(N(x)) \ge p-2$.  Let $H\subseteq N[x]$ be isomorphic to  $K_{p-1}$.  Since $\delta(G)=2p-5$, every vertex in $H$ is adjacent to  $p-3$ vertices in $G\less H$. Then $G\less H$ is disconnected, otherwise $G>K_t$ by contracting $G\less H$ into a single vertex, a contradiction. Let $G_1$ be a component of $G\less H$. Then $N(G_1) :=\{v\in V(H):  v \, \text{ is adjacent to a vertex in}\, G_1\}$ is a minimal separating set of $G$, in particular,  $N(G_1)$ is a clique, contrary to Lemma~\ref{dirac}(ii).  This proves that  $\omega(N(x)) \le p-3$. By \pr{e:falpha},  we see that $\delta(N(x))\ge p-3$.\hfill\vrule height3pt width6pt depth2pt \medskip

\noindent \refstepcounter{counter}\label{e:fchi} (\arabic{counter})
$\chi(N(x)) = p - 2$.

\pf  Suppose to the contrary that $\chi(N(x)) \ne p - 2$.  By \pr{e:falpha}, it is clear that $\chi(N(x)) \ge p - 2$. Thus  $\chi(N(x)) = t$ for some $t \ge p - 1$.
Let $V_1, \dots, V_t$ be the color classes of any proper $t$-coloring of $N(x)$. We may assume that the color classes are ordered so that $V_i = \{a_i\}$ for $i=1,2,\dots,  2t-2p+5$ and $V_j = \{a_j, b_j\}$ for $j=2t-2p+6, \dots, t$. Let $r=2t-2p+6\ge4$. Notice that $\{a_1, a_2, \dots, a_{r-1}\}$ induces a clique in $N(x)$, and so  $r \le p - 2$ by \pr{e:fomega}.
By \pr{e:falpha},  $a_i$ is adjacent to either  $a_j$ or $b_j$ for  each $i\in \{1,2,\dots, r-1\}$ and each $j\in\{r, \dots, t\}$.
Notice that if $t=p-1$, then $r=4$. 
Suppose that $t\ge p$ or that $t=p-1$ and   $a_1, a_2, a_3$   have a common neighbor in $N(x)\less \{a_1, a_2, a_3\}$, say $a_4$. By Lemma~\ref{wonderful} applied to $N(x)$ with $S=\{a_r, b_r\}$ and $M=\{a_{r+1}b_{r+1}, \dots, a_tb_t\}$, there exist $t-r$ pairwise vertex-disjoint paths $P_{r+1}, \dots, P_t$ such that each $P_j$ has ends $\{a_j, b_j\}$ and all its internal vertices in $G \backslash N[x]$. 
By contracting each $P_j$ to a single vertex for all $j\in\{r+1, \dots, t\}$, together with $x, a_1, \dots, a_{r-1}$ if $t \ge p$; and together with $x, a_1, a_2, a_3, a_4$ if $t=p-1$ where $a_4 = a_r$ is a common neighbor of $a_1, a_2, a_3$, we obtain a clique minor on $(t - r) + r = t \ge p$ vertices in the former case and $(t - r) + r + 1 = t +1 = p$ vertices in the latter case,  a contradiction. Thus $t=p-1$ and $a_1, a_2, a_3$   have no common neighbor in $N(x)\less \{a_1, a_2, a_3\}$.

\begin{figure}[htb]
\centering
\includegraphics[width=275px]{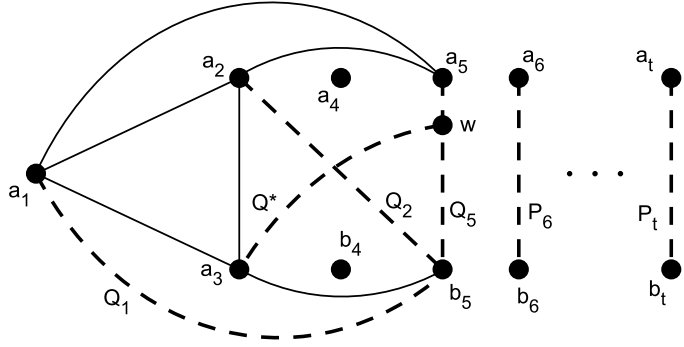}
\caption{Finding a $K_p$ minor when $\chi(N(x)) \ge p - 1$.}
\label{fig:finala1a2a3}
\end{figure}

Since each of $a_1, a_2, a_3$ is adjacent to either $a_5$ or $b_5$, by symmetry, we may assume that  $a_5$ is adjacent to  $a_1$ and $a_2$, but not to $a_3$. Then $b_5$ is adjacent to $a_3$.  We may assume that $b_5$ is not adjacent to $a_1$. For the worst case scenario, we may further assume that $b_5$ is  not adjacent to $a_2$.
By Lemma~\ref{wonderful} applied twice to $N(x)$ with $S=\{a_4, b_4\}$ and $M\in\{ \{\{a_6b_6\}, \dots, \{a_tb_t\}, \{b_5a_1, b_5a_2, b_5a_5\}\}, \{a_5a_3\}\}$,  there exist pairwise vertex-disjoint paths $P_{6}, \dots, P_t$  such that each $P_j$ has ends $\{a_j, b_j\}$ and all its internal vertices in $G \backslash N[x]$; and paths $Q_1, Q_2, Q_5, Q$   with ends $\{b_5, a_1\}$, $\{b_5, a_2\}$,  $\{b_5, a_5\}$ and $\{a_5, a_3\}$, respectively, and all their internal vertices in $G \backslash N[x]$.
Notice that each $P_j$ is vertex-disjoint from $Q_1, Q_2, Q_5, Q$, and $Q$ is vertex-disjoint from $Q_1, Q_2$ but not necessarily from $Q_5$.  Let $w$ be the first vertex on $Q$ (when read from $a_3$ to $a_5$) that is also on $Q_5$.
Note that $w$ could be $a_5$.
Let $Q^*$ be the subpath of $Q$ between $w$ and $a_3$, as depicted in Figure~\ref{fig:finala1a2a3}.
By contracting each $P_j$ to a single vertex for all $j\in\{6, \dots, t\}$,  $Q_1\less a_1$ and $Q_2\less a_2$ onto $b_5$, $Q^*\less w$ onto $a_3$, and $Q_5\less b_5$ onto $a_5$, together with the vertices $x, a_1, a_2, a_3$, we obtain a $K_p$ minor,   a contradiction.   \hfill\vrule height3pt width6pt depth2pt\medskip

\noindent \refstepcounter{counter}\label{e:fdelta} (\arabic{counter})
 $\delta(N(x)) \ge p - 2$.

\pf  By \pr{e:falpha} and \pr{e:fomega}, 
we have $\delta(N(x)) \ge p - 3$. 
Suppose  there exists a vertex $y\in N(x)$ such that $y$ has exactly $p - 3$ neighbors in $N(x)$.
Then $N(x)\less N[y]$ is a clique of $N(x)$ with  $p - 3$  vertices. Furthermore,  by \pr{e:fomega},  $N(x) \cap N(y)$ must have some missing edge, say $uv$.
Then every vertex of $N(x)\less N[y]$ is adjacent to at least one of $u$ and $v$.
Thus, by contracting the path $uyv$ to a single vertex, we can see that $N(x) > K_{p - 2} \cup K_1$.

Now we can assume without loss of generality that $y \in N(x)$ is such that $N(x) \backslash y > K_{p - 2}$.
Clearly $y$ is not adjacent to every vertex in $N(x) \backslash y$, or else $G > N[x] > K_p$, a contradiction.
Let $\{y_1, \dots, y_t\} = N(x) \backslash N[y]$, where $t\ge1$.
Again, by \pr{e:fomega}, $N(y)$ must have some missing edge, say $uv$.
By Lemma~\ref{wonderful} applied  to $N(x)$ with $S = \{u, v\}$ and $M =\{y y_1, \dots, yy_t\}$, there exist  paths $P_1, \dots, P_t$ such that each $P_i$ has ends $\{y, y_i\}$ and all internal vertices in $G \backslash N[x]$.
Now by contracting each $P_i \backslash y_i$ onto $y$, we see that $G > K_p$, a contradiction.
\hfill\vrule height3pt width6pt depth2pt\\

By \pr{e:fomega} and  \pr{e:fdelta}, $N(x)$ does not contain $K_{p - 2}$ as a subgraph and  $\delta(N(x)) \ge p - 2$.
By \pr{e:fchi}, $\chi(N(x))=p-2$.
Let $V_1, \dots, V_{p - 2}$ be the color classes of any proper $(p - 2)$-coloring of $N(x)$. We may assume that the color classes are ordered so that $V_1 = \{a_1\}$ and $V_j = \{a_j, b_j\}$ for $j=2, \dots, p-2$.   By \pr{e:fdelta},  $a_1$ must be  complete to some  color class $V_i\in\{V_2, \dots, V_{p - 2}\}$, say $V_2$.    By  \pr{e:fdelta} again, 
 $a_2$ and $b_2$  must have   one common neighbor in some color class $V_i\in\{V_3, \dots, V_{p - 2}\}$, say  $V_3$. We may further assume that $a_3$ is adjacent to both $a_2, b_2$. By symmetry, we may  assume that $b_3$ is adjacent to $a_2$.  By Lemma~\ref{wonderful} applied up to $N(x)$ with $S = \{a_2, b_2\}$ and $M =\{\{b_3a_1, b_3a_3\}, \{a_4b_4\}, \dots, \{a_{p-2}b_{p-2}\} \}$, 
  there exist paths $P_1, P_2$ and pairwise vertex-disjoint paths $Q_4, \dots, Q_{p -2} $ such that  $P_1, P_2$ have ends $\{b_3, a_1\}$ and $\{b_3, a_3\}$,  respectively;  each $Q_j$ has  ends $\{a_j, b_j\}$ and all such paths have their internal vertices in $G \backslash N[x]$. By contracting  $P_1\less a_1, P_2\less a_3$ onto $b_3$, the edge $b_2a_3$ onto $a_3$, and  each $Q_j$ into a single vertex for $4 \le j \le p - 2$,  we see that $G > K_p$, a contradiction.\medskip
  
This  completes the proof of Theorem~\ref{final}. \hfill\vrule height3pt width6pt depth2pt \\

\noindent {\bf Acknowledgements}

We would like to thank Robin Thomas and Guoli Ding for their helpful discussion. We thank one  referee  for many helpful comments.


\baselineskip 2pt
\begin{thebibliography}{99}


\bibitem{AG2015}
B. Albar,  D. Gon\c calves, 
On triangles in $K_r$-minor free graphs, 	arXiv:1304.5468. 
%
\vspace {-0.25cm}
%
\bibitem{CHS} G. Chen, Z. Hu,  F. Song, A new connectivity bound for linkages and its application to the Hadwiger's  conjecture, submitted.
%
\vspace {-0.25cm}
%
\bibitem{Dirac1952} G. A. Dirac, A property of $4$-chromatic graphs and some 
remarks on critical graphs, J. London Math. Soc. 27 (1952) 85--92.
%
\vspace {-0.25cm}
%
\bibitem{Dirac1964} G. A. Dirac, Homomorphism theorems for graphs, Math. Ann.  153 (1964) 69--80.
%
\vspace {-0.25cm}
%
\bibitem{dirac2}G. A. Dirac, Trennende Knotenpunktmengen und Reduzibilit\"at abstrakter Graphen mit Anwendung auf das Vierfarbenproblem,  J. Reine Agew. Math.  204 (1960) 116--131.
%
\vspace {-0.25cm}
%
\bibitem{hc}H. Hadwiger, \"Uber eine Klassifikation der Strechenkomplexe, 
Vierteljschr. Naturforsch. Ges Z\"urich  88 (1943) 133--142.
%
\vspace {-0.25cm}
%
\bibitem{Jakobsen}
I. T. Jakobsen,
A homomorphism theorem with an application to the conjecture of Hadwiger, Studia Sci. Math. Hungar.  6 (1971) 151--160.
%
\vspace {-0.25cm}
%
\bibitem{Jakobsen1972} I. T. Jakobsen, On certain homomorphism properties of graphs I, 
Math. Scand.  31 (1972) 379--404.
%
\vspace {-0.25cm}
%
\bibitem{Jakobsen1983}I. T. Jakobsen, On certain homomorphism properties of graphs II, 
 Math. Scand.  52 (1983) 229--261.
%
\vspace {-0.25cm}
%
\bibitem{Jorgensen1994}
L. K. J\o rgensen, 
Contractions to $K_8$, 
 J. Graph Theory 18  (1994) 431--448.
%
\vspace {-0.25cm}
%
\bibitem{ktoft}K. Kawarabayashi,  B. Toft, Any $7$-chromatic graph has a $K_7$ or $K_{4,4}$ as a minor, Combinatorica  25 (2005) 327--353.
%
\vspace {-0.25cm}
%
\bibitem{Kawa2007}
K. Kawarabayashi,
On the connectivity of minimum and minimal counterexamples to Hadwiger's Conjecture,
 J. Combin. Theory Ser. B  97 (2007) 144--150.
%
\vspace {-0.25cm}
%
\bibitem{Kawa2013}
K. Kawarabayashi,  G. Yu,
Connectivities for $k$-knitted graphs and for minimal counterexamples to Hadwiger's Conjecture,
J. Combin. Theory Ser. B 103 (2013) 320--326.
%
\vspace {-0.25cm}
%
\bibitem{mader}
W. Mader, 
Homomorphies\"atze f\"ur Graphen,
 Math. Ann. 178 (1968) 154--168.
%
\vspace {-0.25cm}
%
\bibitem{7con}
W. Mader.
\"Uber trennende Eckenmengen in homomorphiekritischen Graphen.
 Math. Ann. 175 (1968)  243--252.
%
\vspace {-0.25cm}
%
\bibitem{RST}N. Robertson, P. Seymour and R. Thomas, 
Hadwiger's conjecture for $K_6$-free graphs, 
{\sl Combinatorica \bf 13} (1993) 279--361.
%
\vspace {-0.25cm}
%
\bibitem{RolekSong2} M. Rolek,  Z-X. Song, Properties of $8$-contraction-critical graphs with no $K_7$ minor, in preparation.
%
\vspace {-0.25cm}
%
\bibitem{Seymour}  P.  Seymour, Hadwiger's conjecture, manuscript.
%
\vspace {-0.25cm}
%
\bibitem{Song2005}
Z.-X. Song,
The extremal function for $K_8^-$ minors,
J. Combin. Theory Ser. B  95 (2005) 300--317.
%
\vspace {-0.25cm}
%
\bibitem{SongThomas2006}
Z-X. Song,  R. Thomas, 
The extremal function for $K_9$ minors.
J. Combin. Theory, Ser.  B, 
96 (2006) 240--252.
%
\vspace {-0.25cm}
%
\bibitem{Toft} B. Toft, {\it A survey on Hadwiger's Conjecture}, in : {\it Surveys in Graph Theory}  (edited by G. Chartrand and M. Jacobson), Congr. Numer. 115 (1996), 249--283.
%
\vspace {-0.25cm}
%
\bibitem{wagner}K. Wagner,  \"Uber eine Eigenschaft der ebenen Komplexe, 
 Math. Ann.  114 (1937) 570--590.





\end{thebibliography}
\end{document}